\documentclass{article}
\usepackage{amsfonts}
\usepackage{amssymb}
\usepackage{amsmath}
\usepackage{amscd}
\usepackage[dvips]{graphicx}
\usepackage{amssymb,amsfonts,amsxtra, mathrsfs,placeins,graphicx,verbatim}
\usepackage[all]{xy}

\newtheorem{theorem}{Theorem}

\newtheorem{definition}[theorem]{Definition}

\newtheorem{lemma}[theorem]{Lemma}
\newtheorem{proposition}[theorem]{Proposition}

\newcommand\C{\mathbb{C}}

\newcommand\g{\mathfrak{g}}

\newcommand\q{\mathcal{A}}

\newcommand\p{\mathcal{P}}

\begin{document}

\title{Poisson structures on $\mathbb{C[}X_{1},\ldots ,X_{n}]$ associated with rigid Lie algebras.}
\author{Nicolas Goze\\
   LMIA,\\
   Universit\'e de Haute Alsace,\\
  4, rue des Fr\`eres Lumi\`ere,\\
  MULHOUSE, 68093. France\\
   \texttt{nicolas.goze@uha.fr}}
\date{\today}

\maketitle

\begin{abstract}
We present the classical Poisson-Lichnerowicz cohomology for the Poisson algebra of polynomials $\mathbb{C}[X_{1},\ldots ,X_{n}]$
using exterior calculus. After presenting some non homogeneous Poisson brackets on this algebra, we compute Poisson cohomological spaces when the Poisson structure corresponds to a bracket of a rigid Lie algebra.
\end{abstract}

\section{Introduction}

The first Poisson structures appeared in classical mechanics. In 1809, D. Poisson
introduced a bracket of functions, which permits to write the Hamilton's
equations as differential equations. This leaded to define a Poisson manifold, that is, a
manifold $M$ whose algebra of smooth functions $F(M)$ is equipped with a skew-symmetric
bilinear map
$$\{\ ,\ \} : F(M) × F(M) \rightarrow F(M),$$
satisfying the Leibniz rule,
$$\{FG, H\} = F\{G, H\} + \{F, H\}G, $$ and the Jacobi identity. In \cite{Lich}, A. Lichnerowicz has also introduced a cohomology, associated to a Poisson structure, called Poisson cohomology.

In this paper we study in terms of exterior calculus the Poisson structures on the associative algebra of complex polynomials in $n$ variables. We apply this approach to the determination of non homogeneous quadratic Poisson brackets and to the computation of the Poisson cohomology. The linear Poisson structures are naturally related to the $n$-dimensional Lie algebras. Recall that a complex Lie algebra $\frak g$ is rigid when its orbit in the algebraic variety of $n$-dimensional complex Lie algebra defined by the Jacobi relations is Zariski open. Such an algebra admits a non trivial Malcev torus and it is  graded by the roots of the torus. We study the Poisson structure on $\mathbb{C}[X_{1},\ldots ,X_{n}]$ whose Poisson brackets correspond to a solvable rigid Lie bracket with non zero roots. In a generic example we compute the corresponding Poisson cohomology.

\section{Poisson structures on $\mathbb{C}[X_{1},\ldots ,X_{n}]$ and
exterior calculus}

\subsection{Poisson brackets and differential forms}

\bigskip Let $\mathcal{A}^{n}$ be the commutative associative algebra $%
\mathbb{C}[X_{1},\ldots ,X_{n}]$ of complex polynomials in $X_{1},\cdots
,X_{n}$. We define a Poisson structure on $\mathcal{A}^{n}$ as a bivector
\begin{equation*}
\mathcal{P}=\sum_{1\leq i<j\leq n}P_{ij}\partial _{i}\wedge \partial _{j},
\end{equation*}%
where $\partial _{i}=\frac{\partial }{\partial X_{i}}$ and $P_{ij}\in
\mathcal{A}^{n}$, satisfying the axiom
\begin{equation*}
\lbrack \mathcal{P},\mathcal{P}]_{S}=0,
\end{equation*}%
where $[,]_{S}$ denotes the Schouten's bracket. If $\mathcal{P}$ is a
Poisson structure on $\mathcal{A}^{n}$, then
\begin{equation*}
\{P,Q\}=\mathcal{P}(P,Q)
\end{equation*}%
defines a Lie bracket on $\mathcal{A}^{n}$ which satisfies the Leibniz identity
\begin{equation*}
\{PQ,R\}=P\{Q,R\}+Q\{P,R\},
\end{equation*}%
for any $P,Q,R\in \mathcal{A}^{n}.$

\noindent We denote by $Sh_{p,q}$ the set of $(p,q)$-shuffles where a $(p,q)$%
-shuffle is a permutation $\sigma $ of the symmetric group $\Sigma_{p+q}$ of
degree $p+q$ such that $\sigma (1)<\sigma (2)<\cdots <\sigma (p)$ and $%
\sigma (p+1)<\sigma (p+2)<\cdots <\sigma (p+q)$. For any bivector $\mathcal{P%
}$ we consider the $(n-2)$-exterior form
\begin{equation*}
\Omega =\sum_{\sigma \in S_{2,n-2}}(-1)^{\varepsilon (\sigma )}P_{\sigma
(1)\sigma (2)}dX_{\sigma (3)}\wedge \cdots \wedge dX_{\sigma (n)},
\end{equation*}%
where $(-1)^{\varepsilon (\sigma )}$ is the signature of the permutation $\sigma $.
If $n>3$, we consider the Pfaffian form $\alpha _{i_{1},\cdots ,i_{n-3}}$
given by
\begin{equation*}
\begin{array}{l}
\alpha _{i_{1},\cdots ,i_{n-3}}(Y)=\Omega (\partial _{i_{1}},\partial
_{i_{2}},\ldots ,\partial _{i_{n-3}},Y)%
\end{array}%
\end{equation*}%
with $Y=\sum_{i=1}^{n}Y_{i}\partial _{i},$ $Y_{i}\in \mathcal{A}^{n}.$

\begin{theorem}
A bivector $\mathcal{P}$ on $\mathcal{A}^{n}$ satisfies $[\mathcal{P},%
\mathcal{P}]_{S}=0$ if and only if

\begin{itemize}
\item for $n>3$,
\begin{equation*}
d\alpha _{i_{1},\cdots ,i_{n-3}}\wedge \Omega =0,
\end{equation*}%
for every ${i_{1},\cdots ,i_{n-3}}$ such that $1\leq i_1 < \cdots <
i_{n-3}\leq n$.

\item for $n=3$
\begin{equation*}
d\Omega\wedge \Omega =0.
\end{equation*}
\end{itemize}
\end{theorem}

\noindent\textit{Proof.} The integrability condition $[\mathcal{P},\mathcal{P%
}]_{S}=0$ writes%
\begin{equation*}
\sum_{r=1}^{n}P_{ri}\partial _{r}P_{jk}+P_{rj}\partial
_{r}P_{ki}+P_{rk}\partial _{r}P_{ij}=0,
\end{equation*}%
for any $1\leq i,j,k\leq n.$ But%
\begin{equation*}
\alpha _{i_{1},\cdots ,i_{n-3}}=\sum (-1)^{N}P_{jk}dX_{l},
\end{equation*}%
summing over all triples $(j,k,l)$ such that $(j,k,i_{1},\ldots
,l,\ldots i_{n-3})$ is a permutation of $S_{2,n-2}$ and $N=\varepsilon
(\sigma )+p-3$ where $(-1)^\varepsilon(\sigma )$ is the signum of $\sigma$. Then
\begin{equation*}
d\alpha _{i_{1},\cdots ,i_{n-3}}=\sum (-1)^{N}dP_{jk}\wedge dX_{l}
\end{equation*}%
and $d\alpha _{i_{1},\cdots ,i_{n-3}}\wedge \Omega =0$ corresponds to $[%
\mathcal{P},\mathcal{P}]_{S}=0.$ The proof is similar if $n=3$.

\subsection{Lichnerowicz-Poisson cohomology}

We denote by ${\mathcal{A}}^n_{\mathcal{P}}$ the algebra $\mathcal{A}^n=%
\mathbb{C[}X_{1},\ldots ,X_{n}]$ provided with the Poisson structure $%
\mathcal{P}. $ For $k \geq 1$, let $\chi ^{k}({\mathcal{A}}^n_{\mathcal{P}})$ be the vector
space of $k$-derivations that is of $k$-skew linear maps on ${\mathcal{A}}%
^n_{\mathcal{P}}$ satisfying%
\begin{equation*}
\varphi (P_{1}Q_{1},P_{2},\ldots ,P_{k})=P_{1}\varphi (Q_{1},P_{2},\ldots
,P_{k})+Q_{1}\varphi (P_{1},P_{2},\ldots ,P_{k}),
\end{equation*}%
for all $Q_{1},P_{1},\ldots ,P_{k}\in {\mathcal{A}}^n_{\mathcal{P}}.$ For $%
k=0$ we put $\chi ^{0}({\mathcal{A}}^n_{\mathcal{P}})={\mathcal{A}}^n_{%
\mathcal{P}}.$ Let $\delta ^{k}$ be the linear map%
\begin{equation*}
\delta ^{k}:\chi ^{k}({\mathcal{A}}^n_{\mathcal{P}})\longrightarrow \chi
^{k+1}({\mathcal{A}}^n_{\mathcal{P}})
\end{equation*}%
given by%
\begin{equation*}
\begin{array}{ll}
\delta ^{k}\varphi (P_{1},P_{2},\ldots ,P_{k+1})=& \displaystyle%
\sum_{i=1}^{k+1}(-1)^{i-1}\{P_{i},\varphi (P_{1},\ldots ,\widehat{P_{i}}%
,\ldots P_{k+1})\} \\
\medskip
&+\displaystyle\sum_{1\leq i<j\leq k+1}(-1)^{i+j}\varphi
(\{P_{i},P_{j}\},P_{1},\ldots ,\widehat{P_{i}},\ldots ,\widehat{P_{j}}%
,\ldots P_{k+1}),
\end{array}%
\end{equation*}%
where $\widehat{P_{i}}$ means that the term $P_{i}$ does not appear. We have
$\delta ^{k+1}\circ \delta ^{k}=0$ and the Lichnerowicz-Poisson cohomology
corresponds to the complex $(\chi ^{k}(\mathcal{A_P}),\delta ^{k})_{k}$. Let
us note that $\chi ^{k}({\mathcal{A}}^n_{\mathcal{P}})$ is trivial as soon
as $k>n.$ A description of the cocycle $\delta ^{k}\varphi $ is presented in
\cite{Pi} for $n=3$ using the vector calculus. We will describe these
formulae using exterior calculus for $n>3$. Let us begin with some notations:

\begin{itemize}
\item To any element $P\in {\mathcal{A}}^n_{\mathcal{P}}=\chi ^{0}({\mathcal{%
A}}^n_{\mathcal{P}}),$ we associate the $n$-exterior form
\begin{equation*}
\Phi _{n}(P)=PdX_{1}\wedge \ldots \wedge dX_{n}.
\end{equation*}

\item To any $\varphi \in \chi ^{k}({\mathcal{A}}_{\mathcal{P}}^{n})$ for $%
1\leq k<n,$ we associate the $(n-k)$-exterior form%
\begin{equation*}
\Phi _{n-k}(\varphi )=\sum_{\sigma \in S_{k,n-k}}(-1)^{\varepsilon (\sigma
)}\varphi (X_{\sigma (1)},\ldots ,X_{\sigma (k)})dX_{\sigma (k+1)}\wedge
\ldots \wedge dX_{\sigma (n)}.
\end{equation*}

\item To any $\varphi \in \chi ^{n}({\mathcal{A}}_{\mathcal{P}}^{n}),$ we
associate the function $\Phi _{0}(\varphi )=\varphi .$
\end{itemize}

\noindent Finally, if $\theta $ is an $k$-exterior form and $%
Y=\displaystyle\sum_{i=1}^{n}Y_{i}\partial _{i}$ is a vector field with $Y_{i}\in
{\mathcal{A}}_{\mathcal{P}}^{n}$, then the inner product $i(Y)\theta $ is the
$(k-1)$-exterior form given by
\begin{equation*}
i(Y)\theta (Z_{1},\cdots ,Z_{k-1})=\theta (Y,Z_{1},\cdots ,Z_{k-1}),
\end{equation*}%
for every vector fields $Z_{1},\cdots ,Z_{k-1}$.

\begin{theorem}
Assume that $n=3$.\ Then we have

1. For all $P\in {\mathcal{A}}_{\mathcal{P}}^{3},$
\begin{equation*}
\Phi _{2}(\delta ^{0}P)=-\Omega \wedge dP.
\end{equation*}

2. For all $f\in \chi ^{1}({\mathcal{A}}_{\mathcal{P}}^{3}),$%
\begin{equation*}
\begin{array}{ll}
\Phi _{1}(\delta ^{1}f)  = & -i(\partial _{1},\partial _{2})[\Omega \wedge
d(i(\partial _{3})\Phi _{2}(f))+d(i(\partial _{3})\Omega )\wedge \Phi
_{2}(f)] \\
& +i(\partial _{1},\partial _{3})[\Omega \wedge d(i(\partial _{2})\Phi
_{2}(f))+d(i(\partial _{2})\Omega )\wedge \Phi _{2}(f)] \\
& -i(\partial _{2},\partial _{3})[\Omega \wedge d(i(\partial _{1})\Phi
_{2}(f))+d(i(\partial _{1})\Omega )\wedge \Phi _{2}(f)],%
\end{array}%
\end{equation*}%
where $i(X,Y)$ denotes the composition $i(X)\circ i(Y).$

3. For all $\varphi \in \chi ^{2}({\mathcal{A}}^3_{\mathcal{P}}),$%
\begin{equation*}
\Phi _{0}(\delta ^{2}\varphi
)=i(\partial_{1},\partial_{2},\partial_{3})(d\Omega \wedge \Phi _{1}(\varphi
)+\Omega \wedge d\Phi _{1}(\varphi )).
\end{equation*}
\end{theorem}

\noindent \textit{Proof.} If $n=3$ we have
\begin{equation*}
\Omega =P_{12}dX_{3}-P_{13}dX_{2}+P_{23}dX_{1}.
\end{equation*}%
Then the integrability of $\mathcal{P}$ is equivalent to $\Omega \wedge
d\Omega =0.$ The theorem results of a direct computation and of the
following general formula:
\begin{equation*}
\forall \varphi \in \chi ^{k}({\mathcal{A}}_{\mathcal{P}}^{n}),\ \varphi
(P_{1},\ldots ,P_{k})=\sum_{1\leq i_{1}\leq \ldots \leq i_{k}\leq n}\varphi
(X _{i_{1}},\ldots ,X _{i_{k}})\partial _{i_{1}}P_{1}\ldots
\partial _{i_{k}}P_{k}.
\end{equation*}

\noindent \textbf{Example.} We consider the Poisson algebra ${\mathcal{A}}_{%
\mathcal{P}_{1}}^{n}=(\mathbb{C[}X_{1},X_{2},X_{3}],\mathcal{P}_{1})$ where $%
\mathcal{P}_{1}$ is given by
\begin{equation*}
\left\{
\begin{array}{l}
\mathcal{P}_{1}(X_{1},X_{2})=X_{2}, \\
\mathcal{P}_{1}(X_{1},X_{3})=2X_{3}, \\
\mathcal{P}_{1}(X_{2},X_{3})=0.%
\end{array}%
\right.
\end{equation*}%
Then
\begin{equation*}
\dim H^{0}({\mathcal{A}}_{\mathcal{P}_{1}}^{n})=1,\ \dim H^{1}({\mathcal{A}}_{%
\mathcal{P}_{1}}^{n})=3,\ \\dim  H^{2}({\mathcal{A}}_{\mathcal{P}_{1}}^{n})=2,\
H^{3}({\mathcal{A}}_{\mathcal{P}_{1}}^{n})=\{0\}.
\end{equation*}%
In this case $\Omega =X_{2}dX_{3}-2X_{3}dX_{2}$ and $d\Omega =3dX_{2}\wedge
dX_{3}$. Let us compute $\dim H^{2}({\mathcal{A}}_{\mathcal{P}_{1}}^{n})$ .
Let $\varphi \in \chi ^{2}({{\mathcal{A}}_{\mathcal{P}}^{n}}_{1}).$ Then\ $%
\Phi _{0}(\delta ^{2}\varphi )=0$ implies%
\begin{equation*}
d\Omega \wedge \Phi _{1}(\varphi )+\Omega \wedge d\Phi _{1}(\varphi ))=0,
\end{equation*}
that is
\begin{equation*}
\begin{array}{l}
X_{2}(\partial _{1}\varphi (X_{1},X_{3})+\partial _{2}\varphi
(X_{2},X_{3}))+2X_{3}(-\partial _{1}\varphi (X_{1},X_{2})+\partial
_{3}\varphi (X_{2},X_{3})) \\
+3\varphi (X_{2},X_{3})=0.%
\end{array}%
\end{equation*}%
Now, if $f\in \chi ^{1}({\mathcal{A}}_{\mathcal{P}_{1}}^{n})$ then
\begin{equation*}
\begin{array}{ll}
\Phi _{1}(\delta f)= & [X_{2}(-\partial _{2}f(X_{2})-\partial
_{1}f(X_{1}))-2X_{3}(\partial _{3}f(X_{2}))+f(X_{2})]dX_{3} \\
& -[2X_{3}(\partial _{1}f(X_{1})+\partial _{3}f(X_{3}))+X_{2}(\partial
_{2}f(X_{3}))-2f(X_{3})]dX_{2} \\
& -[X_{2}(-\partial _{1}f(X_{3}))-2X_{3}(\partial _{1}f(X_{2}))]dX_{1}.%
\end{array}%
\end{equation*}%
Comparing these two relations we obtain that $H^{2}({\mathcal{A}}_{\mathcal{P%
}_{1}}^{n})$ is generated by the two cocycles
\begin{equation*}
\left\{
\begin{array}{l}
\Phi _{1}(\varphi _{1})=X_{3}dX_{2}, \\
\Phi _{1}(\varphi _{2})=X_{2}^{2}dX_{2}. \\
\end{array}%
\right.
\end{equation*}%
\bigskip

Now consider the general case. Let $\mathcal{A}=\mathbb{C[}%
X_{1},\ldots ,X_{n}]$ be provided with the Poisson structure $\mathcal{P}.$

\begin{theorem}
Let $\varphi \in \chi ^{k}(\mathcal{A_{P}})$. Then, we have
$$
\begin{array}{ll}
\Phi _{n-k-1}(\delta ^{k}\varphi )=&\varepsilon \sum i(\partial _{\sigma
(1)},\cdots ,\partial _{\sigma (k+1)})[d(i(\partial _{\sigma (k+2)},\cdots
,\partial _{\sigma (n)})\Omega )\wedge \Phi _{n-k}(\varphi )\\
\medskip\\
&+\Omega \wedge d(i(\partial _{\sigma (k+2)},\cdots ,\partial _{\sigma
(n)})\Phi _{n-k}(\varphi ))],
\end{array}$$
for all $\sigma \in S_{k+1,n-k-1}$, where $\varepsilon =\varepsilon
(n,k)=(-1)^{\frac{(n-k)(n-k+1)}{2}}$.
\end{theorem}
\noindent
\textit{Proof.} To simplify we write $d_{i}$ in place of $dX_{i}$. We have
seen that for every $P\in \mathcal{A_{P}}$ we have $\delta ^{0}P=-\Omega
\wedge dP.$ But
\begin{equation*}
\Phi _{n-1}(\delta P)=\displaystyle\sum_{k=1}^{n}(-1)^{k-1}\{X_{k},P\}d_{1}\wedge
\cdots \wedge \hat{d_{k}}\wedge \cdots \wedge d_{n},
\end{equation*}%
where $\hat{d_{i}}$ means that this factor does not appear, with $\{P,X_{i}\}=%
\displaystyle\sum_{j=1}^{n}P_{ji}\partial _{j}P$ with $P_{ji}=-P_{ij}$ when $j>i$. But
$$
\begin{array}{l}
\smallskip
i(\partial_1)[\Omega \wedge
d(i(\partial_2,\cdots,\partial_{n})\Phi
_{n}(P))+d(i(\partial_2,\cdots,\partial_{n})\Omega )\wedge \Phi _{n}(P)
\\
\smallskip
\ \ \ =i(\partial_1)[\Omega \wedge
d(i(\partial_2,\cdots,\partial_{n})\Phi_{n}(P)] =(-1)^{\frac{n(n-1)}{2}%
}i(\partial_1)[\Omega \wedge dP \wedge d_1]\\

\smallskip
\ \ \ = -(-1)^{\frac{n(n-1)}{2}} \displaystyle\sum_{i=2}^nP_{1i}%
\partial_iPd_2\wedge\cdots\wedge d_n =(-1)^{\frac{n(n-1)}{2}%
}\Phi_{n-1}(P)(\partial_2,\cdots,\partial_n).
\end{array}$$

\noindent Similarly
$$
\begin{array}{l}
\smallskip

 i(\partial_j)[\Omega \wedge d(i(\partial_1,\cdots,\hat{\partial_j}%
,\cdots , \partial_{n})\Phi _{n}(P))+d(i(\partial_1,\cdots,\hat{\partial_j}%
,\cdots ,\partial_{n})\Omega )\wedge \Phi _{n}(P)]
\\
\smallskip

\ \ \ = i(\partial_j)[\Omega \wedge d(i(\partial_1,\cdots,\hat{%
\partial_j},\cdots , \partial_{n})\Phi _{n}(P))=(-1)^{j-1+\frac{n(n-1)}{2}%
}i(\partial_j)[\Omega \wedge dP \wedge dX_j]\\

\smallskip

\ \ \ = (-1)^{j-1+\frac{n(n-1)}{2}}i(\partial_j)(\displaystyle%
\sum_{l=1}^{l=j-1}P_{1j}\partial_lP-\displaystyle\sum_{l=j+1}^{l=n}P_{jl}%
\partial_lP)d_1\wedge\cdots \wedge d_n \\

\smallskip

\ \ \ = (-1)^{\frac{n(n-1)}{2}}(\sum_{l=1}^{l=j-1}P_{1j}\partial_lP-%
\displaystyle\sum_{l=j+1}^{l=n}P_{jl}\partial_lP)d_1\wedge\cdots \wedge \hat{%
d_j}\cdots\wedge d_n \\

\smallskip

\ \ \ =(-1)^{\frac{n(n-1)}{2}}\{P,X_i\}d_1\wedge \cdots\wedge \hat{%
d_i}\wedge \cdots \wedge d_n.

\end{array}$$

\noindent We deduce
\begin{equation*}
\Phi _{n-1}(\delta ^{0}P)=(-1)^{\frac{n(n-1)}{2}}\sum_{j=1}^{n}(-1)^{j-1}i(%
\partial _{j})[\Omega \wedge d(i(\partial _{1},\cdots ,\hat{\partial _{j}}%
,\cdots ,\partial _{n})\Phi _{n}(P))
\end{equation*}%
which proves the theorem for $k=0$. The proof is similar for any $k$.

\bigskip

\bigskip

\noindent \textbf{Application.} We consider the $n$-dimensional complex Lie
algebra defined by the brackets
\begin{equation*}
\lbrack X_{1},X_{i}]=(i-1)X_{i},
\end{equation*}%
for $i=2,\cdots ,n.$ Let $\mathcal{P}_{2}$ be the corresponding Poisson bracket
on $\mathbb{C}[X_{1},\cdots ,X_{n}]$. Let $\chi _{2}^{k}(\mathcal{A}_{%
\mathcal{P}_{2}})$ be the subspace of $\chi ^{k}(\mathcal{A}_{\mathcal{P}%
_{2}})$ whose elements are homogeneous of degree $2$. We denote by $%
H_{2}^{2}(\mathcal{A}_{\mathcal{P}_{2}})=Z_{2}^{2}/B_{2}^{2}$ the
corresponding subspace of $H^{2}(\mathcal{A}_{\mathcal{P}_{2}})$. Define $N:=%
\frac{n(n-1)}{2}$.

\begin{itemize}
\item If $n$ is even, then
\begin{equation*}
\dim  B^2_2=N+(N-1)+\cdots+N-n/2+1=\displaystyle \frac{n(2n^2-3n+2)}{8}.
\end{equation*}

\item If $n$ is odd,
\begin{equation*}
\dim  B^2_2=N+(N-1)+\cdots+(N-(n-1)/2)= \displaystyle \frac{(n^2-1)(2n-1)}{8}.
\end{equation*}
\end{itemize}

\noindent In fact, if $f \in \chi^1_2(\mathcal{A}_{\mathcal{P}})$, then $%
f(X_i)=P_i=\Sigma a_{i_{1},\cdots,i_{n}}^{i}X_{1}^{i_{1}}X_{2}^{i_{2}}\cdots
X_{n}^{i_{n}} $ is homogeneous of degree $2$, then:

1. In $\delta f(X_1,X_{2l})$ we find $N-l$ independent coefficients of $%
P_{2l}$. The coefficients which do not appear are:
\begin{equation*}
a^{2l}_{1,0,0,\cdots,0,1,0,\cdots,0}, \ \
a^{2l}_{0,1,0,\cdots,0,1,0,\cdots,0},\cdots,
a^{2l}_{0,0,\cdots,1,1,0,\cdots,0},
\end{equation*}
where the second $1$ is respectively in the place $2l,2l-1,\cdots,l+1$.

2. In $\delta f(X_1,X_{2l+1})$ we find $N-l-1$ independent coefficients of $%
P_{2l+1}$. The coefficients which do not appear are:
\begin{equation*}
a^{2l+1}_{1,0,0,\cdots,0,1,0,\cdots,0}, \ \
a^{2l+1}_{0,1,0,\cdots,0,1,0,\cdots,0}, \ \ \cdots,
a^{2l+1}_{0,0,..,0,2,0,\cdots,0},
\end{equation*}
where the second $1$ is in place $2l+1,2l,\cdots,l+2$ and in the last case
the $2$ is in place $l+1$.

3. For $i \geq 2$ and $j>i$, $\delta f(X_i,X_j)$ is defined by the $(n-2)$
coefficients $a^i_{1,0,0,\cdots,0,1,0,\cdots,0}$.

\medskip

Now we can to find the generators of $H_{2}^{2}(\mathcal{A}_{%
\mathcal{P}})$. We can choose $\phi \in \chi _{2}^{2}$ such that
\begin{equation*}
\left\{
\begin{array}{l}
\phi (X_{1},X_{2})=0, \\
\phi (X_{1},X_{3})=a_{1,3}^{1,3}X_{1}X_{3}+a_{1,3}^{2,2}X_{2}^{2}, \\
\cdots \\
\phi
(X_{1},X_{2l})=a_{1,2l}^{1,2l}X_{1}X_{2l}+a_{1,2l}^{2,2l-1}X_{2}X_{2l-1}+%
\cdots +a_{1,2l}^{l,l+1}X_{l}X_{l+1}, \\
\phi
(X_{1},X_{2l+1})=a_{1,2l+1}^{1,2l+1}X_{1}X_{2l+1}+a_{1,2l+1}^{2,2l}X_{2}X_{2l}+\cdots +a_{1,2l+1}^{l,l}X_{l}^{2},
\\
\cdots \\
\phi
(X_{1},X_{n})=a_{1,n}^{1,n}X_{1}X_{n}+a_{1,n}^{2,n-1}X_{2}X_{n-1}+\cdots ,
\\
\phi (X_{i},X_{j})=A_{i,j},%
\end{array}%
\right.
\end{equation*}%
where $A_{i,j}$ is a degree $2$ homogeneous polynomial without monomial of
type $X_{1}X_{k}$ and $X_{i}X_{j}$. By solving $\Phi _{n-2}(\delta \phi )=0$
we obtain the generators of $H_{2}^{2}(\mathcal{A}_{\mathcal{P}})$. They are
given by
\begin{equation*}
\left\{
\begin{array}{l}
\phi (X_{1},X_{2})=0, \\
\phi (X_{1},X_{3})=a_{1,3}^{2,2}X_{2}^{2}, \\
\cdots \\
\phi (X_{1},X_{2l})=a_{1,2l}^{2,2l-1}X_{2}X_{2l-1}+\cdots
+a_{1,2l}^{l,l+1}X_{l}X_{l+1}, \\
\phi (X_{1},X_{2l+1})=a_{1,2l+1}^{2,2l}X_{2}X_{2l}+\cdots
+a_{1,2l+1}^{l+1,l+1}X_{l+1}^{2}, \\
\cdots \\
\phi (X_{1},X_{n})=a_{1,n}^{2,n-1}X_{2}X_{n-1}+\cdots
+a_{1,n}^{m,m+1}X_{m}X_{m+1},\ \ \mbox{\rm if}\ n=2m, \\
\phi (X_{i},X_{j})=A_{i,j},%
\end{array}%
\right.
\end{equation*}%
\noindent or $\phi (X_{1},X_{n})=a_{1,n}^{2,n-1}X_{2}X_{n-1}+\cdots
+a_{1,n}^{m+1,m+1}X_{m}X_{m+1},\ \ \mbox{\rm if}\ n=2m+1.$ \medskip
\noindent For example:

- if $n=2$ , $\dim  H^2_2( \mathcal{A}_{\mathcal{P}_2}, \mathcal{A}_{\mathcal{P%
}_2})=1$,

- if $n=3$, $\dim  H^2_2( \mathcal{A}_{\mathcal{P}_2}, \mathcal{A}_{\mathcal{P}%
_2})=3$,

- if $n=4$, $\dim  H^2_2( \mathcal{A}_{\mathcal{P}_2}, \mathcal{A}_{\mathcal{P}%
_2})=8$,

- if $n=5$, $\dim H_{2}^{2}(\mathcal{A}_{\mathcal{P}_{2}},\mathcal{A}_{%
\mathcal{P}_{2}})=16$.

\section{Poisson structures of degree $2$ on $\mathbb{C}[X_1,X_2,X_3]$}

\bigskip Let $\mathcal{P} $ be a Poisson structure on $\mathcal{A}^3=\mathbb{C}%
[X_{1},X_{2},X_{3}]$ with $P_{ij}$ of degree $2$. Then $\mathcal{P}$  writes
\[
\mathcal{P}=\mathcal{P}_{0}+\mathcal{P}_{1}+\mathcal{P}_{2},
\]%
where $\mathcal{P}_{i}$ is homogeneous of degree $i$. The associated form $%
\Omega $ is decomposed in homogeneous parts $\Omega =\Omega _{0}+\Omega
_{1}+\Omega _{2}$ and, since $d\Omega_0=0$, the condition $\Omega \wedge d\Omega=0$ is equivalent to
\begin{equation}
\left\{
\begin{array}{l}
\medskip \Omega _{2}\wedge d\Omega _{2}=0, \\
\medskip \Omega _{0}\wedge d\Omega _{1}+\Omega _{1}\wedge d\Omega _{0}=0, \\
\medskip \Omega _{0}\wedge d\Omega _{2}+\Omega _{2}\wedge d\Omega
_{0}+\Omega _{1}\wedge d\Omega _{1}=0, \\
\medskip \Omega _{1}\wedge d\Omega _{2}+\Omega _{2}\wedge d\Omega _{1}=0.%
\end{array}
\right.%
\end{equation}%
If $\Omega _2=0$, then $\p$ is a linear Poisson structure on $\q ^3$ (\cite{JPD}). If $\Omega _2 \neq 0$ and
$\Omega_0=\Omega _1=0$, then $\p$ is a quadratic homogeneous Poisson structure and the classification is given in \cite{Mo}.
In this section we will study the remaining cases
$\Omega_0 \neq 0$ or $\Omega _1\neq 0$. The associative algebra $\q ^3$ admits a natural grading
$\q ^3= \oplus _{n\geq 0}V_n$ where $V_n$ is the space of degree $n$ homogeneous polynomial of $\q ^3$.

\begin{definition}
A linear isomorphism
\[
f:\oplus _{n\geq 0}V_{n}\rightarrow \oplus _{n\geq 0}V_{n}
\]%
is called equivalence of order $2$ if it satisfies
\begin{itemize}
\item $f(V_{1})\subset V_{1}\oplus V_{2},$

\item $f(V_{0})=V_{0},$

\item $f\mid _{\oplus _{n\geq 2}V_{n}}=Id.$
\end{itemize}
Moreover if $V_1$ is provided with a Lie algebra structure, then
\begin{itemize}
\item
$\pi_{1}\circ f$ is a Lie automorphism of $V_{1}$,
\end{itemize}
where $\pi _{1}$ is the projection on $V_{1}.$
\end{definition}
Such a map  writes%
\[
\left\{
\begin{array}{l}
f(X_{i}) =\displaystyle \sum_{j=1}^n a_{i}^{j}X_{j}+\sum_{j,k=1}^n b_{i}^{jk}X_{j}X_{k}, \\
f(X_{i}X_{j}) =X_{i}X_{j}.
\end{array}
\right.%
\]
Thus, if $\p$ is a degree $2$ Poisson structure on $\q ^3$, putting $Y_{i}=f(X_{i})$ and
\[
\{Y_{i},Y_{j}\}=f^{-1}(\{f(X_{i}),f(X_{j})\}),
\]%
we obtain a new Poisson structure of degree $2$. These two Poisson structures are called equivalent. In the following, we classify
the non homogeneous Poisson structure of degree $2$ up to an equivalence of order $2$. Note that the quadratic homogeneous Poisson structures are classified in \cite{Ha}. We assume also that these Poisson structures are not trivial extensions of Poisson structures on $\mathcal{A}^2$, that is, Poisson structures which do not depend only of two variables.

\subsection{First case: $\Omega = \Omega _2 + \Omega _1, \Omega _1 \neq 0$}
The integrability condition of $\Omega$ reduces to
\begin{equation}
\left\{
\begin{array}{l}
\medskip \Omega _{1}\wedge d\Omega _{1}=0,\\
\medskip \Omega _{1}\wedge d\Omega _{2}+\Omega _{2}\wedge d\Omega _{1} =0,\\
\medskip \Omega _{2}\wedge d\Omega _{2} =0.
\end{array}%
\right.%
\end{equation}
As $\Omega ^{1}\wedge d\Omega ^{1}=0$, $\Omega _1$ defines on $\q ^3$ a linear Poisson structure. Then this form is isomorphic to
one of the following
\[
\left\{
\begin{array}{l}
\Omega ^{1}_1=X_{3}dX_{3}, \\
\Omega ^{2}_1=X_{2}dX_{3}+X_{3}dX_{2}+X_{1}dX_{1}, \\
\Omega ^{3}_1=X_{2}dX_{3}-\alpha X_{3}dX_{2}, \\
\Omega ^{4}_1=(X_{2}+X_{3})dX_{3}-X_{3}dX_{2}. \\
\end{array}
\right.%
\]
Consider $
\Omega _{2}=A_{3}dX_{3}-A_{2}dX_{2}+A_{3}dX_{1}
$  with%
$$
\left\{
\begin{array}{lll}
A_{1}
&=&a_{1}X_{1}^{2}+a_{2}X_{2}^{2}+a_{3}X_{3}^{2}+a_{4}X_{1}X_{2}+a_{5}X_{1}X_{3}+a_{6}X_{2}X_{3},
\\
A_{2}
&=&b_{1}X_{1}^{2}+b_{2}X_{2}^{2}+b_{3}X_{3}^{2}+b_{4}X_{1}X_{2}+b_{5}X_{1}X_{3}+b_{6}X_{2}X_{3},
\\
A_{3}
&=&c_{1}X_{1}^{2}+c_{2}X_{2}^{2}+c_{3}X_{3}^{2}+c_{4}X_{1}X_{2}+c_{5}X_{1}X_{3}+c_{6}X_{2}X_{3}.
\end{array}
\right.
$$
\subsubsection{ $d\Omega _{2}=0$}
If $ \Omega _1 = \Omega _1 ^1$ or $\Omega _1^2$, then $d\Omega _1 =0$ and (3.2) is satisfied. An equivalence of order $2$
of type $Y_1=X_1,Y_2=X_2,Y_3=X_3+B$ where B is an homogeneous polynomial of degree $2$, allows to reduce the form $\Omega_2$ to a form with $A_1=0$. We obtain the following
Poisson structure associated to
\begin{equation}
\begin{array}{l}
\Omega(1)= (aX_{1}^{2} -\frac{b}{2}X_{2}^{2}-2cX_{1}X_{2})dX_1
-(cX_{1}^{2}+eX_{2}^{2}+bX_{1}X_{2})dX_2+X_3dX_3
\end{array}
\end{equation}
corresponding to $ \Omega _1 = \Omega _1 ^1$, and
\begin{equation}
\begin{array}{ll}
\Omega(2)= &(X_1+aX_{1}^{2}-\frac{b}{2}X_{2}^{2}-2cX_{1}X_{2})dX_1+(X_3-cX_{1}^{2}-eX_{2}^{2} -bX_{1}X_{2})dX_2\\
& + X_3dX_3
\end{array}
\end{equation}
corresponding to $ \Omega _1 = \Omega _1 ^2$. If $ \Omega _1 = \Omega _1 ^3$ or $\Omega _1^4$, then
$d\Omega _1 =kdX_2\wedge dX_3$ with $k \neq 0$. Then (3.2) implies $\Omega_2 \wedge dX_2\wedge dX_3=0$ that is
$A_3=0$. Such a structure is a Poisson structure on $\mathcal{A}^2$.

\subsubsection{ $d\Omega _{2}\neq 0$, $\Omega _1=\Omega_1^1$}
As $d\Omega _{1}= 0$, then (3.2) is equivalent to
\[
\left\{
\begin{array}{l}
\Omega _{1}\wedge d\Omega _{2}=0, \\
\Omega _{2}\wedge d\Omega _{2}=0.
\end{array}%
\right.
\]%
This implies $Pd\Omega _{2}=\Omega _{1}\wedge \Omega _{2}$ where $P$ is an
homogeneous polynomial of degree $2$.
The equivalence of order $2$ given by $%
Y_{i}=X_{i} $ for $i=1,2$ and $Y_{3}=X_{3}+B$ with $B\in V_{2}$ enables to consider $A_1=0$.
In this case, $\Omega _{1}\wedge \Omega
_{2}=Pd\Omega _{2}$ is equivalent to
\[
\left\{
\begin{array}{l}
\partial _{1}A_{2}+\partial _{2}A_{3}=0, \\
P\partial _{3}A_{3}=X_{3}A_{3}, \\
P\partial _{3}A_{2}=X_{3}A_{2}.%
\end{array}%
\right.
\]%
If $X_3$ is not a factor of $P$, then $\partial_{3}A_{2}=\alpha X_{3}$ and $%
\partial _{3}A_{3}=\beta X_{3}.$ If $\alpha =\beta =0,$ then $\Omega_2=0$.
The case $\alpha \beta \neq 0$ reduces by a change of variables to the case
 $\alpha \neq 0$ and $\beta =0$, then $A_3=0$. Thus we obtain
$$
\Omega= X_3dX_3 - (aX_{2}^{2}+bX_{3}^{2})dX_2.
$$
This structure is a trivial extension of a Poisson structure on $\mathbb{C}[X_2,X_3]$.
If $P=X_{3}Q$ and $Q$ is a degree $1$ homogeneous
polynomial, then $Q$ satisfies
\[
\left\{
\begin{array}{l}
Q(\partial _{1}A_{2}+\partial _{2}A_{3})=0, \\
Q\partial _{3}A_{3}=A_{3}, \\
Q\partial _{3}A_{2}=A_{2}.%
\end{array}%
\right.
\]%
We deduce the following structures
$$
\begin{array}{l}
\Omega=  (aX_1^2+bX_1X_3)dX_1 +X_3dX_3, \\
\Omega=  (aX_1+X_3/2)^2dX_1+X_3dX_3, \\
\Omega=  (aX_1X_3+bX_2X_3)dX_1+(bX_1X_3+cX_2X_3)dX_2+X_3dX_3,\\
\end{array}
$$
The two first ones depends only of two variables. Then we obtain the following Poisson structure
\begin{equation}
\Omega(3)=  (aX_1X_3+bX_2X_3)dX_1+(bX_1X_3+cX_2X_3)dX_2+X_3dX_3.
\end{equation}

\subsubsection{ $d\Omega _{2}\neq 0$, $\Omega _1=\Omega_1^2$}
By an equivalence of degree $2$, we can consider that $A_3=0$. Then $Pd\Omega_2=\Omega_1 \wedge \Omega_2$ gives
$$
\left\{
\begin{array}{l}
P\partial _{1}A_{2}=X_{1}A_{2}, \\
P\partial _{1}A_{1}=X_{1}A_{1}, \\
P(\partial _{2}A_{1}+\partial _{3}A_{2})=(A_{2}X_{2}+A_{1}X_{3}).
\end{array}
\right.
$$
Solving these equations, we obtain:
\begin{equation}
\begin{array}{l}
\Omega(4)= X_1dX_1+(X_3-aX_1X_3)dX_2+(X_2+aX_1X_2)dX_3,\\
\Omega(5)=  X_1dX_1 +(X_3-aX_1^2-2aX_2X_3)dX_2+X_2dX_3.
\end{array}
\end{equation}

\subsubsection{$d\Omega _{2}\neq 0$, $\Omega _1=\Omega_1^3=X_{2}dX_{3}-\alpha X_{3}dX_{2}$}
Assume that $\alpha \neq 0$ and $\alpha \neq -1$. The
equivalence given by $Y_{2}=X_{2}+B_{2},Y_{i}=X_{i}$ for $i=1,3$ and $%
B_{2}\in V_{2}$ shows that the structure corresponding to $\Omega=\Omega_1$
is equivalent to a structure of degree $2$ defines as follow
\[
\left\{
\begin{array}{l}
A_{1}=a_{2}X_{2}^{2}+a_{3}X_{3}^{2}+\frac{c_{6}}{\alpha }%
X_{1}X_{2}+c_{3}X_{1}X_{3}, \\
A_2=0, \\
A_{3}=c_{3}X_{3}^{2}+c_{5}X_{1}X_{3}+c_{6}X_{2}X_{3}.
\end{array}
\right.
\]
Thus we can assume that in $\Omega _2$ we have $%
c_{3}=c_{5}=c_{6}=a_{2}=a_{3}=a_6=0$. The new equivalence of degree $2$
given by $Y_{3}=X_{3}+B_{3},Y_{i}=X_{i}$ for $i=1,2$ and $B_{3}\in V_{2}$
gives a Poisson structure of degree $2$ equivalent to the structure of
degree $1$ with
\[
\left\{
\begin{array}{l}
A_1=0, \\
A_{2}=b_{2}X_{2}^{2}+b_{3}X_{3}^{2}-c_{2}X_{1}X_{2}+\frac{c_{6}}{\alpha }%
X_{1}X_{3}, \\
A_{3}=c_{2}X_{2}^{2}+c_{4}X_{1}X_{2}+c_{6}X_{2}X_{3}.
\end{array}
\right.
\]
Thus we can assume that
$$
\Omega_2=(a_{1}X_{1}^{2}+a_{4}X_{1}X_{2}+a_{5}X_{1}X_{3})dX_1+(b_{1}X_{1}^{2}+b_{4}X_{1}X_{2}+b_{5}X_{1}X_{3})dX_2+c_{1}X_{1}^{2}dX_3.
$$
 As $\Omega_{1}\wedge d\Omega _{2}+\Omega _{2}\wedge d\Omega _{1}=0$ we obtain the following Poisson structure
\begin{equation}
\begin{array}{l}
\Omega(6)=aX_1X_3dX_1-\alpha X_3dX_2+(X_{2}-\frac{a}{2\alpha}X_1^2)dX_3.
\end{array}
\end{equation}
with $\alpha \neq 0$ and $\alpha \neq -1$.

\medskip

\noindent If $\alpha =-1$, then $d\Omega _1 =0$ and this case has already been studied. If $\alpha =0,$ \ by equivalence of degree $2$ we can assume that
\[
\left\{
\begin{array}{l}
A_1=a_1X_1^2+a_3X_3^2+a_4X_1X_2+a_5X_1X_3, \\
A_{2}=b_{1}X_{1}^{2}+b_{3}X_{3}^{2}+b_5X_1X_3, \\
A_{3}=c_{1}X_{1}^{2}+c_{3}X_3^2+c_{5}X_{1}X_{3}.
\end{array}
\right.
\]
Then we have $A_3=0$ and the Poisson structure concerns only two variables.

\subsubsection{$d\Omega_2 \neq 0,  \Omega_1=\Omega_1^4=(X_{2}+X_{3})dX_{3}-X_{3}dX_{2}$}

 By
equivalence of degree $2$, we can assume that $A_{1}=0,c_{5}=0$ and $%
b_{4}=0. $ The equation $\Omega _{1}\wedge d\Omega_{2}+\Omega _{2}\wedge
d\Omega _{1}=0$ implies that $c_{1}=c_{4}=b_{1}=0,c_{6}=-b_{5}=2c_{2}.$ The
equation $\Omega _{2}\wedge d\Omega _{2}=0$ implies that $c_{2}=0$ and $%
b_{2}c_{3}=b_{6}c_{3}=0.$ Then we obtain the following Poisson structure:

\noindent\begin{equation}
\Omega(7)=  aX_3^2dX_1 - (X_3+bX_3^2)dX_2+(X_2+X_3)dX_3,  \\
\end{equation}
with $a \neq 0$.
\subsection{Second case: $\Omega = \Omega _2 + \Omega _1 +\Omega _0 , \Omega _0 \neq 0$}
The form $\Omega_0 \oplus \Omega_1$ provides the vector space $V_{0}\oplus
V_{1}$ with a linear Poisson structure. Then $V_{0}\oplus V_{1}$ is a Lie
algebra such that $V_0$ is in the center. This implies $\Omega _{1}\wedge
d\Omega _{1}=0.$ We deduce that  $\Omega_0 + \Omega_1$ is equivalent to
\begin{equation}
\left\{
\begin{array}{l}
\medskip dX_{3}-X_{3}dX_{2}, \\
\medskip X_{3}dX_{3}-dX_{2}, \\
\medskip X_{2}dX_{3}+X_{3}dX_{2}+dX_{1}.
\end{array}
\right.%
\end{equation}

\subsubsection{ $\Omega_0 + \Omega _1=dX_{3}-X_{3}dX_{2} $}
By equivalence, we can assume that $%
a_{3}=a_{5}=b_{3}=c_{5}=0.$ The equation $\Omega _{0}\wedge d\Omega _{2}=0$
implies $b_{4}=-2c_{2},c_{4}=-2b_{1},c_{6}=-b_{5}$, $\Omega _{1}\wedge
d\Omega _{2}+\Omega _{2}\wedge d\Omega _{1}=0 $ implies that $
c_{1}=c_{2}=c_{3}=c_{4}=0,a_{1}=a_{4}=0$ and $\Omega _{2}\wedge d\Omega
_{2}=0$ gives $b_{5}b_{2}=b_{5}a_{2}=b_{5}a_{6}=0.$ Thus we obtain the following Poisson structures given by
\noindent\begin{equation}
\Omega(8)= aX_2X_3dX_1-(X_3-aX_1X_3+bX_2X_3)dX_2+dX_3  \\
\end{equation}
with $a \neq 0$.
\subsubsection{ $\Omega_0 + \Omega _1=-dX_{2}+X_{3}dX_{3} $}

We can assume that $%
A_{2}=b_{2}X_{2}^{2}+b_{4}X_{1}X_{2}.$ As $d\Omega _{1}=0$, the system
reduces to $\Omega _{0}\wedge d\Omega _{2}=\Omega _{1}\wedge d\Omega _{2}=0$%
. This gives $c_{4}=c_{6}=a_{4}=0$ and $%
b_{4}+2c_{2}=a_{5}-2c_{3}=2a_{1}-c_{5}=0.$ Thus $\Omega _{2}\wedge d\Omega_{2}=0$
is equivalent to $(2a_{2}X_{2}+a_{6}X_{3})A_{3}=0.$ We obtain the following Poisson structures
\noindent\begin{equation}
\begin{array}{ll}
\Omega(9)=   &- (X_2+aX_2^2+bX_1X_2)dX_2+(1+cX_1^2+eX_3^2+fX_1X_3)dX_3 \\ &+(gX_1^2-\frac{b}{2}X_2^2
+\frac{f}{2}X_3^2+2cX_1X_3)dX_1. \\
\end{array}
\end{equation}

\subsubsection{$\Omega_0+\Omega _1=dX_1+X_{3}dX_{2}+X_{2}dX_{3}$}
 By equivalence, we can assume $%
b_{5}=b_{2}=a_{3}=a_{5}=c_{2}=c_{5}=0.$ As $d\Omega _{1}=0$, the equation $%
\Omega _{0}\wedge d\Omega _{2}=\Omega _{1}\wedge d\Omega _{2}=0$ implies that $%
b_{6}+2a_{2}=a_{6}+2b_{3}=a_{4}=a_{1}=b_{1}=b_{4}=c_{3}=0 $ . In this case $%
\Omega _{2}\wedge d\Omega _{2}=0$ is equivalent to $%
c_{6}(X_{2}A_{2}+X_{3}A_{1})=0.$ We obtain
\noindent\begin{equation}
\Omega(10)=   (1+aX_1^2)dX_1+X_3dX_2+X_2dX_3,
\end{equation}
and
\begin{equation}
\Omega(11)=  (1+aX_1^2)dX_1+(X_3+bX_3^2+cX_2X_3)dX_2+(X_2+\frac{c}{2}X_2^2+2bX_2X_3)dX_3.
\end{equation}

\section{Poisson algebras associated to rigid Lie algebras}

\subsection{Rigid Lie algebras}

\bigskip Let us fix a basis of $\mathbb{C}^{n}.$ With respect to this basis, a
multiplication $\mu $ of a $n$-dimensional complex Lie algebra is determined
by its structure constants $C_{ij}^{k}$. We denote by $L_{n}$ the algebraic variety \ $\mathbb{C}[C_{ij}^{k}]/I$
where $I$ is the ideal generated by the polynomials:%
\[
\left\{
\begin{array}{l}
C_{ij}^{k}+C_{ji}^{k}=0, \\
\displaystyle\sum_{l=1}^{n}C_{ij}^{l}C_{lk}^{s}+C_{jk}^{l}C_{li}^{s}+C_{ki}^{l}C_{li}^{s}=0,
\end{array}
\right.
\]%
for all $1\leq i,j,k,s\leq n.$ Then every multiplication $\mu $ of a $n$%
-dimensional complex Lie algebra is identified to one point of $L_{n}.$ We
have a natural action of the algebraic group $Gl(n,\mathbb{C)}$ on $L_{n}$
whose orbits correspond to the classes of isomorphic multiplications:%
\[
\mathcal{O(\mu )=}\left\{ f^{-1}\circ \mu \circ (f\times f),\quad f\in Gl(n,%
\mathbb{C)}\right\} .
\]
Let $\mathfrak{g}=(\mathbb{C }^n,\mu) $ be a $n$-dimensional complex Lie
algebra. We denote also by $\mu $ the corresponding point of $L_n$.

\medskip

\begin{definition}
The Lie algebra $\mathfrak{g}$ is rigid if its orbit $\mathcal{O}(\mu
)$ is open (for the Zariski topology) in $L_n$.
\end{definition}

\medskip

Among rigid complex Lie algebras, there are all simple and semi-simple
Lie algebras, all Borel algebras and parabolic Lie algebras. Concerning the classification of rigid Lie algebras, we know
the classification up the dimension $8$ (\cite{G.A}), the classification in any dimension  of solvable rigid Lie algebras whose nilradical is
filiform (\cite{G.A}). Recall two interesting tools to study rigidity of a given Lie algebra.

\medskip

\begin{theorem}
Let $\mathfrak{g}=(\mathbb{C}^n,\mu)$ be a $n$-dimensional complex Lie algebra.
Then

1. $\mathfrak{g}$ is rigid if and only if any valued deformation $\mathfrak{g}^{\prime }$ is
($K^*$)-isomorphic to $\mathfrak{g}$ where $K^*$ is the fraction field of the valuation ring $R$ containing the structure constants of
$\mathfrak{g}^{\prime }$.

2. (Nijenhuis-Richardson Theorem) If $H^2(\mathfrak{g},\mathfrak{g})=0$, then $\mathfrak{g}$ is rigid.
\end{theorem}

The notion of valued deformation, which extends in a natural way the classical
notion of Gerstenhaber deformations, is
developed in \cite{G.R1}. In the Nijenhuis-Richardson theorem, the second cohomological space $H^2(\mathfrak{g},\mathfrak{g})$ of the Chevalley
cohomology of $\g$ is trivial. Let us recall that the converse of this  theorem is not true. There exists solvable rigid Lie algebras with $H^{2}(%
\mathfrak{g,g})\neq 0$ (see for example \cite{G.A}). In this case there exists a $2$-cocycle $\varphi_1 \in  H^2(\mathfrak{g},\mathfrak{g})$ which is not
the first term of a valued (or formal) deformation \[
\mu _{t}=\mu +\sum_{i\geq 1}t^{i}\varphi _{i}
\]
of the Lie multiplication $\mu$ of $\g$.

\subsection{Finite dimensional Poisson algebras whose Lie bracket is rigid}

We recall in this section some results of \cite{G.R2} which precise the structure of a finite dimensional complex Poisson algebra
with rigid underlying Lie bracket.
Let $\p=(\C^n,\p)$ be a finite dimensional complex Poisson algebra. We denote by $\{X,Y\}$ and $X \cdot Y$ the corresponding Lie bracket and
associative multiplication, by $\g_{\p}$ the Lie algebra $(\p,\{,\})$ and by $\mathcal{A}_{\p}$ the associative
algebra $(\p,\cdot)$.
\begin{proposition}(\cite{G.R2}).
If the Lie algebra $\g_{\p}$ is a simple complex Lie algebra, then the associative product is trivial that is  $X\cdot Y=0$ for every $X$, $Y$ in $\mathcal{P}$.
\end{proposition}
Let us assume now that $\g_{\p}$ is a complex rigid solvable Lie algebra.
Then $\mathfrak{g}$ is written:
\[
\mathfrak{g}=\mathfrak{t} \oplus \mathfrak{n},
\]
where $\mathfrak{n}$ is the nilradical of $\mathfrak{g}$ and $\mathfrak{t}$
a maximal abelian subalgebra such that the adjoint operators $adX$ are
diagonalizable for every $X \in \mathfrak{t}$. This subalgebra $\mathfrak{t}$ is usually called a Malcev torus. All these
maximal torus are conjugated and their common dimension is called the rank
of $\mathfrak{g}$.
\begin{lemma}
\label{lemma2}
If there is a non-zero vector $X \in \frak{g}_{\mathcal{P}}$ such that $adX$ is
diagonalizable
with $0$ as a simple root, then $\mathcal{A}_\mathcal{P}.\mathcal{A}_\mathcal{P}=\{0\}$.
\end{lemma}

\noindent {\it Proof.} Let $\left\{ e_1, ..., e_n \right\}$ be a
basis of $\frak{g}_\mathcal{P}$ such that $ad  e_1$ is diagonal
with respect to this basis. By assumption, $\{e_1,e_i\}= \lambda_i
e_i$ with $\lambda_i \neq 0$ for $i \geq 2$. Since $\{
e_1^2,e_1\}=2e_1 \cdot \{e_1,e_1\}= 0$, it follows that
$e_1^2=ae_1.$ But for any $i \neq 1,$ $\{e_1^2,e_i\}=2e_1 \cdot \{e_1,e_i\}=2 \lambda_i e_1 \cdot e_i$ and
$\{e_1^2,e_i\}=a\lambda_ie_i$, thus $e_1\cdot e_i=\frac{a}{2}e_i$. The associativity of the product $X\cdot Y$
implies that $(e_1 \cdot e_1)\cdot  e_i=ae_1 \cdot
e_i=\frac{a^2}{2}e_i =e_1 \cdot (e_1 \cdot e_i)=\frac{a^2}{4}e_i$. Therefore $a=0$ and $e_1^2=0=e_1 \cdot e_i$ for any $i.$
Finally, $0=\{e_1 \cdot e_j,e_i\}=e_1\cdot
\{ e_j,e_i\}+e_j \cdot \{e_1 ,e_i\}=\lambda_i e_j \cdot e_i$,
which implies $ e_i\cdot e_j=0, \ \forall i,j \geq 1.$ $\Box$
\begin{proposition}
Let $\frak{g}$ be a rigid solvable Lie algebra of rank 1  with non-zero
roots. Then there is only one Poisson algebra $\p$
such that $\frak{g}_\mathcal{P}=\frak{g}$.
It corresponds to $$ X\cdot Y =  0,$$
for any $X,Y \in \p$.
\end{proposition}
\noindent {\it Proof.} By hypothesis we have $\dim  \, \frak{t}=1$ and for $X \in
\mathfrak{g}_\mathcal{P}$, $X \neq 0$, as the roots of $\frak{g}$ are non zero, the restriction of the
operator $ad X$ on $\frak{n}$ is invertible (all known
solvable rigid Lie algebras satisfy this hypothesis). By the previous lemma,
the associated algebra $\mathcal{A}_\mathcal{P}$ satisfies
$\mathcal{A}_\mathcal{P}.\mathcal{A}_\mathcal{P}= \{0\}$.
\begin{theorem}
\label{theo22}
Let $\mathcal{P}$ a complex Poisson algebra such that
$\frak{g}_\mathcal{P}$ is  rigid solvable  of rank 1 (i.e $\dim  \, \frak{t}=1$) with
non-zero roots. Then
$\mathcal{P}$ is a rigid Poisson algebra.
\end{theorem}
\noindent
{\it Proof. } See \cite{G.R2}

\subsection{Linear Poisson structures on $\mathcal{A}^{n+1}=\C[X_0,X_1,\cdots,X_n]$ given by  a rigid Lie bracket}
In this section we consider a linear Poisson bracket on $\C[X_0,\cdots,X_n]$ such that the brackets $\{X_i,X_j\}=\p(X_i,X_j)$ corresponds to a solvable
rigid Lie algebra $\g$ of rank $1$. We assume that the roots (see \cite{G.A}) of this rigid Lie algebras are $1,\cdots,n$. In this case we have
\[
\left\{
\begin{array}{l}
\{X_0,X_i\}=iX_i, \ \ i=1,\cdots,n \\,
\{X_1,X_i\}=X_{i+1}, \ \ i=2,\cdots,n-1 \\,
\{X_2,X_i\}=X_{i+2}, \ \ i=3, \cdots, n-2.
\end{array}
\right.
\]
We denote this $(n+1)$-dimensional Poisson algebra by $\mathcal{P}(\g)$. This algebra is a deformation of the Poisson algebra studied in Section 1.2. The corresponding $(n-1)-$exterior form is
$$
\begin{array}{ll}
\Omega & =\sum\limits_{i=1}^{n}(-1)^{i-1}X_{i}d_{1}\wedge \cdots\wedge \hat{d_{i}}\wedge \cdots\wedge
d_{n}+\sum\limits_{i=2}^{n-1}(-1)^{i}X_{i+1}d_{0}\wedge d_{2}\wedge
\cdots\wedge \hat{d_{i}}\wedge \cdots\wedge
d_{n}\\
&+\sum\limits_{i=3}^{n-2}(-1)^{i+1}X_{i+2}d_{0}\wedge d_{1}\wedge
d_{3}\wedge \cdots\wedge \hat{d_{i}}\wedge \cdots\wedge d_{n},
\end{array}%
$$
where $d_{i}$ denotes $dX_{i}$ and $\hat{d_{i}}$ means that this
term does not appear. Let $\varphi$ be a $2$-cochain. We denote by $\varphi(i,j)$ the vector $\varphi(X_i,X_j)$.
Then $\varphi $ is a $2$ cocycle if and only if%
\begin{equation*}
\begin{array}{lll}
\Phi _{n-1}(\varphi ) & = & (-1)^{n-2}\varphi (1,i)d_{0}\wedge d_{2}\wedge
\cdots\wedge \hat{d_{i}}\wedge \cdots\wedge d_{n} \\
&& +  \sum\limits_{i=3}^{n}(-1)^{i-1}\varphi (2,i)d_{0}\wedge d_{1}\wedge
d_{3}\wedge \cdots\wedge \hat{d_{i}}\wedge \cdots\wedge d_{n} \\
&& +  \sum\limits_{3\leq i<j\leq n}(-1)^{j-i-1}\varphi (i,j)d_{0}\wedge
\cdots\wedge \hat{d_{i}}\wedge \cdots\wedge \hat{d_{j}}\wedge
\cdots\wedge d_{n}%
\end{array}%
\end{equation*}%
satisfies%
\begin{equation}
d[i(\partial _{\sigma (1)},\cdots,\partial _{\sigma (n-2)})\Omega ]\wedge \Phi
_{n-1}(\varphi )+\Omega \wedge d[i(\partial _{\sigma (1)},\cdots,\partial
_{\sigma (n-2)})\Phi _{n-2}(\varphi )]=0,
\end{equation}
for any $\sigma \in S_{3,n-2}.$
As $\g=\frak{t}\oplus \frak{n}$, we have the decomposition $\mathcal{P}(\g)=\mathcal{P}(\frak{t})\oplus \mathcal{P}(\frak{n})$
 where $\mathcal{P}(\frak{t})$ and $\mathcal{P}(\frak{n})$) are the Poisson algebras $(\C[X_0],\p)$
and $(\C[X_{1},\cdots,X_{n}],\p)$. From the Hochschild-Serre
factorization theorem,  we assume that the cocycles are $\frak{t}$-invariant
and with values in $\mathcal{P}(\frak{n})$. We denote this space by $\chi
^{k}(\mathcal{P}(\frak{g}),\mathcal{P}(\frak{g}))^{\frak{t}}.$ If $\ f\in \chi ^{1}(\mathcal{P}(\frak{g}),\mathcal{P}(\frak{g}))^{\frak{t}}$
then $$
\{X_{0},f(X_{i})\}=if(X_{i}),$$
and we obtain
$$f(X_{1})=a_{1}^{1}X_{1},f(X_{2})=a_{1}^{11}X_{1}^{2}+a_{2}^{2}X_{2},\cdots,f(X_{i})=\sum\limits_{l_{1}+\cdots+l_{k}=i}a_{i}^{l_{1}
\cdots l_{k}}X_{1}^{l_{1}}\cdots X_{k}^{l_{k}}.$$
Thus $\delta
f(X_{1},X_{i})=a_{1}^{1}\{X_{1},X_{i}\}+\{X_{1},f(X_{i})\}-f(X_{i+1})$ and we can reduce any element $\varphi \in
Z^{2}(\mathcal{P}(\frak{g}),\mathcal{P}(\frak{g}))^{\frak{t}}$ to a $2$-cocycle satisfying%
\begin{equation*}
\varphi (X_{1},X_{i})=0\text{ for }i=2,\cdots,n-1.
\end{equation*}
We denote by $Z_{k}^{\ast }(\mathcal{P}(\frak{g}),\mathcal{P}(\frak{g}))^{\frak{t}}$  the subspace of homogeneous cocycles of degree $k.$  Let us look the system on the $\varphi(i,j)$ which is deduced from Equation(4.1).

- If ($\sigma (1),\cdots,\sigma (n-2))=(3,4,\cdots,n)$ then Condition (4.1) is trivial.

- If ($\sigma (1),\cdots,\sigma (n-2))=(2,3,\cdots,\hat l,\cdots,n)$ then Condition (4.1) is trivial  as soon as $l \neq n$. If $l=n$ we obtain
$$n\varphi (X_{1},X_{n})+(-1)^{n-1}\sum iX_{i}\partial _{i}\varphi(X_{1},X_{n})=0$$ and $\varphi (X_{1},X_{n})$ is of weight $n+1.$

- If ($\sigma (1),\cdots,\sigma (n-2))=(1,2,\cdots,\hat i,\cdots,\hat j,\cdots,n)$ we obtain
$$(i+j)\varphi (X_{i},X_{j})=\sum kX_{k}\partial
_{k}\varphi (X_{i},X_{j})$$ and $\varphi (X_{i},X_{j})$ is of weight $i+j.$

Other relations show that the space of cocycles of degree $2$ is generated by $\varphi(X_1,X_n)$ and
$\varphi(X_2,X_{2k+1})$ with $k=1,\cdots,l$ where $n=2l+1$ or $n=2l$. The relations between these generators leads to study two cases: $k=1$ and $k=2.$

\bigskip

\paragraph{Case $k=1$.}

As $\varphi (X_{1},X_{n})$ is of weight $n+1$, then $\varphi (X_{1},X_{n})=0.$ We have also  $\varphi (X_{i},X_{j})=a_{ij}^{i+j}X_{i+j}$ if $i+j\leq n.$

If ($\sigma (1),\cdots,\sigma (n-2))=(1,2,\cdots,\hat i,\cdots,\hat j,\cdots,n)$ we obtain
$$(i+j)\varphi (X_{i},X_{j})=\sum kX_{k}\partial
_{k}\varphi (X_{i},X_{j})$$ and $\varphi (X_{i},X_{j})$ is of weight $i+1.$
Then
$$\varphi (X_{i},X_{j})=a_{ij}^{i+j}X_{i+j}$$ if $i+j\leq n.$

If ($\sigma (1),\cdots,\sigma (n-2))=(0, 1,2,\cdots,\hat i,\cdots,\hat j,\cdots,\hat k,\cdots,n)$ with $i\geq 3$, then the related conditions are always satisfied.

If ($\sigma (1),\cdots,\sigma (n-2))=(0,3,\cdots, \hat i,\cdots, n)$, $i \geq 3$,  we obtain
relation between $\varphi (3,l)$ and $\varphi (2,l+1).$ We deduce that
$$a_{3,i}=-a_{2,i+1} + a_{2,i}$$
and  $a_{2,3}=a_{2,4}$.

If ($\sigma (1),\cdots,\sigma (n-2))=(0, 1,2,\cdots,\hat i,\cdots,\hat j,\cdots,\hat k,\cdots,n)$ with $i\geq 3$ then
$$a_{4,i}=a_{2,i+2}-2a_{2,i+1}+a_{2,i}$$
and
$$a_{3,4}=a_{3,5}.$$

If ($\sigma (1),\cdots,\sigma (n-2))=(0,2,3,\cdots, \hat i,\cdots, \hat j,\cdots, n)$, $i \geq 4$, then we have
$$a_{i+1,j}=-a_{i,j+1}+a_{i,j}$$
and
$$a_{i,i+2}=a_{i,i+1}.$$

If ($\sigma (1),\cdots,\sigma (n-2))=(0,1,3,\cdots, \hat i,\cdots, \hat j,\cdots, n)$, $i \geq 4$, then  we have
$$a_{i+2,j}=-a_{i,j+2}+a_{i,j}$$
and
$$a_{3,j}=a_{2,j}-a{2,j+1}.$$
If we solve this linear system, we obtain
\begin{proposition}
If $n\geq 7$, then $H_{1}^{2}(A_{p},A_{p})$ is of dimension $1$ and
generated by the cocycle given by%
\begin{equation*}
\left\{
\begin{array}{l}
\varphi (X_{2},X_{i})=(4-i)X_{2+i}\text{ \ }i=5,\cdots,n-2, \\
\varphi (X_{3},X_{i})=X_{3+i}\text{ \ }i=4,\cdots,n-3, \\
\varphi (X_{i},X_{j})=0\text{ \ in other cases}.%
\end{array}%
\right.
\end{equation*}
\end{proposition}

\paragraph{Case $k=2$.}

The set of generators is of dimension $\frac{p^2+5p}{2}$ if $n=2p+1$ and $\frac{p^2+3p-2}{2}$ if $n=2p$.
The number of independent relations concerning these parameters is greater than the dimension of the set of generators as soon as $n \geq 6$. For $n=5$ the dimension is equal to $2$ and for $n=6$, this dimension is $0$.
We deduce that $\dim  H^2_2=0$ when $n \geq 7$.

\bigskip

\noindent{\bf Remark: Deformations of the Enveloping algebra of a rigid Lie algebra} Let $\mathfrak{g}$ be a finite dimensional complex Lie algebra. We denote by
$\mathcal{U}(\mathfrak{g})$ its enveloping algebra. One of the most
important problem in this time is to look the deformations of the
associative algebra $\mathcal{U}(\mathfrak{g})$. The theory of quantum
groups comes from the deformation of $\mathcal{U}(sl(2))$. In this case, $\mathfrak{g}=sl(2)$ is a rigid Lie algebra and $\mathcal{U}(sl(2))$ is a
rigid associative algebra. Thus we have to look what happens for any rigid
Lie algebra. The aim of this section is to study the deformations of $%
\mathcal{U}(\mathfrak{g})$ when $\mathfrak{g}$ is the rigid Lie algebras
studied in the previous paragraph.

We denote by $S(\mathfrak{g)}$ the symmetric algebra on the vector space $%
\mathfrak{g.}$ This associative commutative algebra is interpreted as the
algebra of polynomials on the dual vector space $\mathfrak{g}^{\ast }$ of $%
\mathfrak{g}$ that is $\mathbb{C}[\alpha _{1},\cdots ,\alpha _{n}]$ where $%
\{\alpha _{1},\cdots ,\alpha _{n}\}$ is a basis of $\mathfrak{g}^{\ast }$.
But the Lie structure of $\mathfrak{g}$ induces a Linear Poisson structure
(or of degree $1$), $\mathcal{P}$, on $\mathfrak{g}^{\ast }.$ In fact, if $%
\{X_{1},\cdots ,X_{n}\}$ is the (dual) basis of $\mathfrak{g}$, this Poisson
structure corresponds to the Poisson structure of degree $1$ on $\mathbb{C}%
[X_{1},\ldots ,X_{n}]$ associated to $\mathfrak{g.}$ From the formality
theorem of Kontsevich, $\mathcal{U}(\mathfrak{g})$ is a deformation of the
Poisson algebra $(\mathbb{C}[X_{1},\ldots ,X_{n}],\mathcal{P)}.$ In his
thesis, Toukaidine Petit (\cite{Pe}) shows that every nontrival deformation of
the Poisson structure $\mathcal{P}$ on $\mathbb{C}[X_{1},\ldots ,X_{n}]$
induces a nontrivial deformation of the associative algebra $\mathcal{U}(%
\mathfrak{g}).$ As a consequence, we have that if $\mathfrak{g}$ is a
nonrigid Lie algebra, then there is a nontrivial deformation of $\mathcal{U}(%
\mathfrak{g}).$

If we consider the rigid Lie algebra $\frak{g}_{n+1}$ studied in the previous paragraph, we have determinate a non trivial cocycle of degree one for the corresponding
Poisson algebra which is not integrable. Thus  we cannot define a deformation of its enveloping algebra. But the Lie algebra $\frak{g}_{n+1}$
admit a deformation in the following nonLie algebra which is written
\begin{equation*}
\left\{
\begin{array}{l}
\mu (X_{0},X_{i})=iX_{i}\text{ \ }i=1,\cdots,n, \\
\mu (X_{1},X_{i})=X_{i+1}\text{ \ }i=2,\cdots,n-1, \\
\mu (X_{2},X_{3})=X_{5}, \\
\mu (X_{2},X_{i})=(5-i)X_{2+i}\text{ \ }i=4,\cdots,n-2, \\
\mu (X_{3},X_{i})=X_{3+i}\text{ \ \ }i=4,\cdots,n-3.%
\end{array}%
\right.
\end{equation*}%

\end{document}